\documentclass{amsart}

\usepackage{geometry,amsmath,amscd,amssymb,enumerate,bbm,latexsym,mathtools}
\usepackage{graphicx}
\usepackage{float}
\usepackage{hyperref}

\usepackage{kpfonts}  
\usepackage[T1]{fontenc}
\usepackage[cp1250]{inputenc}
\usepackage[polish,english]{babel}

\newtheorem*{corollary}{Corollary}
\newtheorem{theorem}{Theorem}
\newtheorem*{claim}{Claim}

\makeatletter
\def\th@plain{%
	\thm@notefont{}
	\itshape 
}
\def\th@definition{%
	\thm@notefont{}
	\normalfont 
}
\def\th@procedure{
	\thm@notefont{}
	\ttfamily
}
\makeatother

\theoremstyle{definition} 
\newtheorem{example}{Example}
\newtheorem{remark}{Remark}

\newcommand{\assign}{:=}

\newcommand{\cdummy}{\,}
\newcommand{\nin}{\not\in}
\newcommand{\nocomma}{}
\newcommand{\tmem}[1]{{\em #1\/}}
\newcommand{\tmop}[1]{\ensuremath{\operatorname{#1}}}

\newenvironment{enumerateroman}{\begin{enumerate}[(i)] }{\end{enumerate}}

\newenvironment{quoteenv}{\begin{quote} }{\end{quote}}

\newenvironment{claimproof}{\noindent\underline{Proof of Claim.}\,}{\leavevmode\unskip\penalty9999 \hbox{}\nobreak\hfill\quad\hbox{$\Diamond$}\medskip}

\author{Szymon Brzostowski}
\author{Tadeusz Krasi\'nski}
\author{Justyna Walewska}
\title[Arnold's problem on monotonicity of the Newton number]{Arnold's problem on monotonicity of the Newton number for surface singularities}
\address{Wydzia\l{} Matematyki i Informatyki, Uniwersytet \L{}\'o{}dzki, Banacha 22, 90-238 \L{}\'o{}d\'z{}, Poland}
\email{brzosts@math.uni.lodz.pl}
\email{krasinsk@uni.lodz.pl}
\email{walewska@math.uni.lodz.pl}

\begin{document}

\maketitle

\begin{abstract}
	
	According to the Kouchnirenko Theorem,
	for a generic (precisely non-degenerate in  the Kouchnirenko sense) isolated singularity $f$ its
	Milnor number $\mu (f)$ is equal to the Newton number $\nu (\Gamma_{+}(f))$ of a combinatorial object associated to $f$, the
	Newton polyhedron $\Gamma_+ (f)$. We give a simple condition characterising,
	in terms of $\Gamma_+ (f)$ and $\Gamma_+ (g)$, the equality $\nu (\Gamma_{+}(f)) = \nu (\Gamma_{+}(g))$, for any
	surface singularities $f$ and $g$ satisfying $\Gamma_+ (f) \subset \Gamma_+
	(g)$. This is a complete solution to an Arnold's problem (1982-16) in this
	case.
	
\end{abstract}

\section{Introduction}
\label{intro}

Let $f : (\mathbbm{C}^n, 0) \rightarrow (\mathbbm{C}, 0)$ be an
{\tmem{isolated singularity}} (that is $f$ possesses an isolated critical
point at $0 \in \mathbbm{C}^n$), in the sequel: a {\tmem{singularity}}, in 
short. The Milnor number $\mu (f)$ (see \cite{Mil68}) of a generic $f$ can be expressed,
as proved by Kouchnirenko \cite{Kou76}, using a combinatorial object associated to $f$,
the Newton polyhedron $\Gamma_+ (f) \subset \mathbbm{R}^n_{\geqslant
	0}$. More precisely, under an appropriate non-degeneracy condition imposed on
$f$, it holds $\mu (f) = \nu (
\Gamma_{+}(f))$, where $\nu (\Gamma_{+}(f))$ is the {\tmem{Newton
		number of $\Gamma_{+}(f)$}}. For $\Gamma_+ (f)$ {\tmem{convenient}} (which means that the Newton polyhedron contains a point on each coordinate axis) the latter number is equal to
\[ \nu (\Gamma_{+}(f)) \assign n! \ V_n - (n - 1) ! \ V_{n - 1} + \ldots + (-
1)^{n - 1} \ 1! \ V_1 + (- 1)^n\ V_0, \]
where $V_n$ is the $n$-dimensional volume of the (usually non-convex)
polyhedron `under' $\Gamma_+ (f)$, $V_{n - 1}$ is the sum of $(n -
1)$-dimensional volumes of the polyhedra `under' $\Gamma_+ (f)$ on all
hyperplanes $\{ x_i = 0 \}$, $V_{n - 2}$ is the sum of $(n - 2)$-dimensional
volumes of the polyhedra `under' $\Gamma_+ (f)$ on all hyperplanes $\{ x_i =
x_j = 0 \}$, and so on.

In his acclaimed list of problems, V.~I.~Arnold posed the following
({\cite[1982-16]{Arn04}}):

\begin{quoteenv}
	`Consider a Newton polyhedron $\Delta$ in $\mathbbm{R}^n$ and the number
	$\mu (\Delta) = n!V - \Sigma (n - 1) !V_i + \Sigma (n - 2) !V_{i \nocomma j}
	- \cdots$, where $V$ is the volume under $\Delta$, $V_i$ is the volume under
	$\Delta$ on the hyperplane $x_i = 0$, $V_{i \nocomma j}$ is the volume under
	$\Delta$ on the hyperplane $x_i = x_j = 0$, and so on. \\
	Then $\mu (\Delta)$ grows (non strictly
	monotonically) as $\Delta$ grows (whenever $\Delta$ remains coconvex
	{\tmem{and integer?}}). {\tmem{There is no elementary proof even for $n =
			2$}}.'
\end{quoteenv}

{\noindent}(here, Arnold's terminology slightly differs from ours: $\Delta$ should
be understood as\label{czy jeszcze inaczej?} $\mathbbm{R}^n_{\geqslant 0}~\setminus~\Gamma_+ (f)$ for a singularity
$f$, and then $\mu (\Delta) = \nu (\Gamma_{+}(f))$). In the comments to the problem,
S.~K. Lando wrote that the monotonicity of $\mu (\Delta)$ follows from the
semi-continuity of the spectrum of a singularity, proved independently by
A.~N.~Varchenko \cite{Var84} and by J.~Steenbrink \cite{Ste85}, and that he himself had given an
elementary proof for $n = 2$ (unpublished). Such a proof (for $n=2$) was eventually
published by A.~Lenarcik \cite{Len08}. In the case of an arbitrary $n$, other
proofs were offered by M.~Furuya \cite{Fur04},
J.~Gwo{\'z}dziewicz \cite{Gwo08} and C.~Bivi{\`a}-Ausina \cite{Biv09}.

In the present paper we essentially complete the solution of the
problem for surface singularities, i.e.~for $n = 3$. More specifically, we
prove not only the monotonicity property of the Newton number, but also we
give a simple geometrical condition characterising the situations in which
this monotonicity is strict. We may describe this condition in the following
intuitive way (for the precise statement see Theorem \ref{Main}): for any $f, g :
(\mathbbm{C}^3, 0) \rightarrow (\mathbbm{C}, 0)$ such that $\Gamma_+ (f)
\subset \Gamma_+ (g)$ one has $\nu (\Gamma_+ (f)) = \nu (\Gamma_+ (g))$ if,
and only if, $\Gamma_+ (f)$ and $\Gamma_+ (g)$ differ by (possibly several)
pyramids with bases in the coordinate planes and heights equal to $1$. The
proof we propose is purely geometrical and elementary. We believe that a
similar result should be valid in the $n$-dimensional case, but, if one simply
tries to mimic the proof offered here, the amount of new combinatorial
complications seems to increase enormously. We also expect that our result (and
its potential multidimensional generalization) will have interesting
applications in many aspects of effective singularity theory, e.g.:
computation of the {\L}ojasiewicz exponent, jumps of the Milnor numbers in
deformations of singularities, searching for tropisms of "partial" gradient ideals
$\left( \frac{\partial f}{\partial z_1}, \ldots, \widehat{\frac{\partial
		f}{\partial z_i}}, \ldots, \frac{\partial f}{\partial z_n} \right)\mathcal{O}_n$ of an
isolated singularity $f$, etc.

Similar problem characterizing $f$ for which $\mu(f)$ is minimal (and equal to $\nu(\Gamma_{+}(f))$) among all singularities with the same Newton polyhedron $\Gamma_{+}(f)$ is given in the recent paper by P.~Mondal \cite{Mon16}.

\section{Polyhedra}
\label{sec:Polyhedra}
According to the standard definitions (see e.g.~M.~Berger \cite{Ber09}), a {\tmem{convex
		$n$-polyhedron in $\mathbbm{R}^n$}} is an intersection of a finite family of
{closed} half-spaces of $\mathbbm{R}^n$, having non-empty interior. An
{\tmem{$n$-polyhedron in $\mathbbm{R}^n$}} is a union of finitely many convex
$n$-polyhedra in $\mathbbm{R}^n$. Let $k \leqslant n$; a {\tmem{$k$-polyhedron
		in $\mathbbm{R}^n$}} is a finite union of $k$-polyhedra in $k$-dimensional
affine subspaces of $\mathbbm{R}^n$. A compact connected $k$-polyhedron in
$\mathbbm{R}^n$ is called a {\tmem{$k$-polytope in $\mathbbm{R}^n$}}.

For convenience, we introduce the following notations. Let $\mathbbm{P},
\mathbbm{Q} \subset \mathbbm{R}^n$ be two $k$-polyhe\-dra. The {\tmem{polyhedral
		difference}} ({\tmem{p-difference}}) of $\mathbbm{P}$ and $\mathbbm{Q}$ is the
closure of their set-theoretical difference, in symbols
\[ \mathbbm{P}-\mathbbm{Q} \assign \overline{\mathbbm{P}\setminus\mathbbm{Q}}
. \]
One can check that $\mathbbm{P}-\mathbbm{Q}$ is also a $k$-polyhedron in
$\mathbbm{R}^n$, or an empty set. We say that $\mathbbm{P}$ and $\mathbbm{Q}$
are {\tmem{relatively disjoint}}, if their relative interiors are disjoint (in
appropriate $k$-dimensional affine subspaces). In particular, two
$n$-polyhedra in $\mathbbm{R}^n$ are relatively disjoint if their interiors
are disjoint.

We define the Newton polyhedra in an abstract way without any relation to singularities. A subset $\Gamma_+ \subset
\mathbbm{R}_{\geqslant 0}^n$ is said to be a {\tmem{Newton polyhedron}} when
there exists a subset $A \subset \mathbbm{N}_0^n$ such that
\[ \Gamma_+ = \tmop{conv} \left( \bigcup_{\mathbbm{i} \in A}
(\mathbbm{i}+\mathbbm{R}_{\geqslant 0}^n) \right) . \]
For such an $A$ we will write $\Gamma_+ = \Gamma_+ (A)$. In the sequel we will
assume that are no superfluous points in $A$, implying $A$ is precisely 
the set of all the vertices of $\Gamma_+$.

\begin{remark}
	In the context of singularity theory, we take $A = \tmop{supp} f$, where $f
	= \sum_{\mathbbm{i} \in \mathbbm{N}^n_0} a_{\mathbbm{i}} z^{\mathbbm{i}}$
	around $0$ and $\tmop{supp} f \assign \{ \mathbbm{i} \in \mathbbm{N}_0^n :
	a_{\mathbbm{i}} \neq 0 \}$.
\end{remark}

A Newton polyhedron $\Gamma_+$ is called {\tmem{convenient}} if $\Gamma_+$ intersects all coordinate axes of
$\mathbbm{R}^n$. Since $\mathbbm{N}_0^n$ is a lattice in
$\mathbbm{R}^n_{\geqslant 0}$, the boundary of a convenient polyhedron
$\Gamma_+$ is a finite union of convex $(n - 1)$-polytopes (compact
{\tmem{facets}}) and a finite union of convex unbounded $(n - 1)$-polyhedra
(unbounded {\tmem{facets}}) lying in coordinate hyperplanes. By $\Gamma$ we
will denote the set of these compact facets, and sometimes -- depending on the
context -- also their set-theoretic union. The closure of the complement of
$\Gamma_+$ in $\mathbbm{R}^n_{\geqslant 0}$ will be denoted by $\Gamma_-$,
i.e.
\[ \Gamma_- \assign \mathbbm{R}^n_{\geqslant 0} - \Gamma_+ . \]
It is an $n$-polytope in $\mathbbm{R}^n$ provided $\Gamma_{-}\neq\varnothing$. Hence, $\Gamma_-$ has finite
$n$-volume. Similarly, for any $\varnothing \neq I \subset \{ 1, \ldots, n
\}$, $\Gamma_-$ restricted to the coordinate hyperplane
$\mathbbm{R}^I_{\geqslant 0} \assign \{ (x_1, \ldots, x_n) \in
\mathbbm{R}^n_{\geqslant 0} : \hspace{0.25em} x_i = 0 \text{ for } i \nin I
\}$, that is $\Gamma_-^I \assign \Gamma_- \cap \mathbbm{R}^I_{\geqslant
	0}$, has finite $(\#I)$-volume. Consequently, we may define the {\tmem{Newton
		number $\nu (\Gamma_+)$ of convenient $\Gamma_+$}} by the formula
\[ \nu (\Gamma_+) \assign n!V_n - (n - 1) ! \ V_{n - 1} + \ldots + (- 1)^{n
	- 1} \ 1! \ V_1 + (- 1)^n\ V_0, \]
where $V_i$ denotes the sum of $i$-volumes of $\Gamma_-^I$, for all
$\varnothing \neq I \subset \{ 1, \ldots, n \}$ satisfying $\#I = i$. Note that $V_0 =1$ if $\Gamma_{-}\neq\varnothing$ and $V_0=0$ if $\Gamma_{-}=\varnothing$. Hence $\nu(\mathbb{R}^n_{\geqslant 0})=0$. Clearly,
we may also extend the domain of this definition to any $n$-polytope
$\mathbbm{P}$ in $\mathbbm{R}^n$; thus $\nu (\mathbbm{P})$ makes sense. Then
for any Newton polyhedron $\Gamma_+$ we have $\nu (\Gamma_+) = \nu
(\Gamma_-)$. We will use both notations interchangeably.

The following notions will be useful in our proof. Let $B$ be a compact $(n - 1)$-polyhedron  in an $(n -
1)$-dimensional hyperplane $H \subset \mathbbm{R}^n$ and $Q \in
\mathbbm{R}^n \setminus H$. A {\tmem{pyramid $\mathbbm{P} (B, Q)$ with apex
		$Q$ and base $B$}} is by definition the cone, $\tmop{cone} (B, Q)$, with
vertex $Q$ and base $B$. 
\begin{figure}[H]
	\centering
	\includegraphics{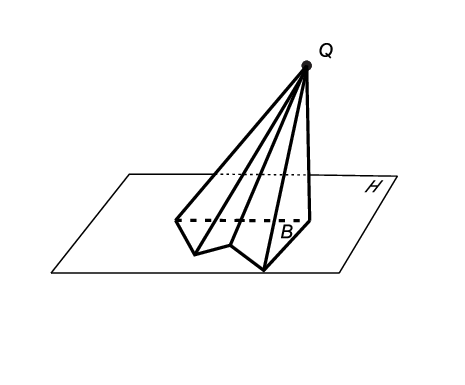}\\
	\caption{The pyramid with apex $Q$ and base $B$.}\label{rys:1}
\end{figure}
\noindent By ({\cite[12.2.2, p.~13]{Ber09}}), the $n$-volume of $\mathbbm{P} (B,
Q)$ can be computed using the elementary formula
\begin{equation}
\tmop{vol}_n \mathbbm{P} (B, Q) = \frac{\tmop{vol}_{n - 1} (B) \cdummy
	\tmop{dist} (Q, H)}{n} . \label{for1}
\end{equation}
Finally, if $P_0, \ldots, P_k \in \mathbbm{R}^n$ are linearly independent points, then
by $\Delta (P_0, \ldots, P_k)$ we denote the {\tmem{$k$-simplex with vertices
		$P_0, \ldots, P_k$}}.

\section{The main Theorem}

Let $\Gamma_+$, $\widetilde{\Gamma}_+$ be two convenient Newton polyhedra such
that $\Gamma_+ \varsubsetneq \widetilde{\Gamma}_+$. Then
\[ \widetilde{\Gamma}_+ = \tmop{conv} (\Gamma_+ \cup \{ P_1, \ldots, P_k \}), \]
for some points $P_1, \ldots, P_k$ lying under $\Gamma_+$, i.e. $P_i \in
\mathbbm{N}_0^n \setminus \Gamma_+$. In such situation $\widetilde{\Gamma}_+$
will also be denoted by $\Gamma_+ + \{ P_1, \ldots, P_k \}$ or $\Gamma_+^{P_1,
	\ldots, P_k}$. Clearly, $\Gamma_+ + \{ P_1, \ldots, P_k \} = (\Gamma_+ + \{
P_1, \ldots, P_{k - 1} \}) + P_k$ and hence $\nu (\Gamma_+ + \{ P_1, \ldots,
P_k \}) = \nu ((\Gamma_+ + \{ P_1, \ldots, P_{k - 1} \}) + P_k)$. Since
moreover
\[ \nu (\Gamma_+) - \nu (\Gamma_+ + \{ P_1, \ldots, P_k \})  = \sum_{1
	\leqslant i \leqslant k} (\nu
(\Gamma_+ + \{ P_1, \ldots, P_{i - 1} \})-\nu (\Gamma_+ + \{ P_1, \ldots, P_i \})), \]
it suffices to consider the monotonicity of the Newton number
for polyhedra defined by sets which differ in one point only, i.e. for Newton
polyhedra $\Gamma_+$ and $\Gamma_+^P$, for some $P \in \mathbbm{N}_0^n
\setminus \Gamma_+$.

\begin{theorem}
	\label{Main}Let $\Gamma_+$ be a convenient Newton polyhedron in
	$\mathbbm{R}^3_{\geqslant 0}$ and let a lattice point $P$ lie under
	$\Gamma_+$ i.e. $P \in \mathbbm{N}_0^3 \setminus \Gamma_+$. Then
	\begin{enumerate}
		\item \label{it1}$\nu (\Gamma_+^P) \leqslant \nu (\Gamma_+)$,
		
		\item \label{it2}$\nu (\Gamma_+^P) = \nu (\Gamma_+)$ if, and only if,
		there exists a coordinate plane $H$ such that $P \in H$ and $\Gamma_+^P -
		\Gamma_+$ is a pyramid with base $(\Gamma_+^P - \Gamma_+) \cap H$ and of
		height equal to $1$.
	\end{enumerate}
\end{theorem}

\begin{figure}[H]
	\centering
	\includegraphics{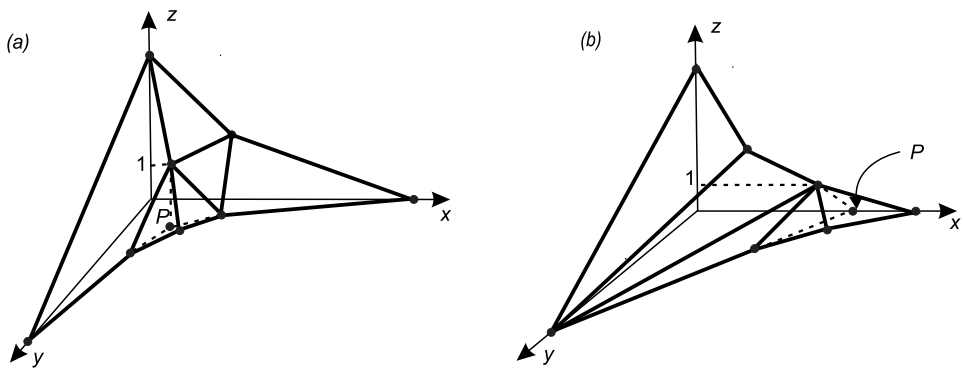}\\
	\caption{$(a)$ $P$ lies in the plane and not on axes.  $(b)$ $P$ lies in the plane and on  an axis.}\label{rys:2}
\end{figure}

\begin{remark}
	We believe that the same theorem is true, {\tmem{mutatis mutandis}}, in the
	$n$-dimensional case. In the simpler case $n = 2$, the theorem is well-known
	(\cite{Gwo08}, \cite{Len08} or \cite{GN12}).
\end{remark}
\begin{remark}
	In the particular case when $\Gamma_+^{P}-\Gamma_+$ is a $3$-dimensional simplex item \ref{it2}  follows from Lemma (2.2) in \cite{Oka89}. 
\end{remark}
\begin{example}
	Let us illustrate the second item of the theorem with some figures. Let $P$
	lying under $\Gamma_+$ be such that $\nu (\Gamma_+^P) = \nu (\Gamma_+)$. Up
	to permutation of the variables, we have the following, essentially
	different, possible locations for $P$:
	\begin{enumerate}
		\item $P$ lies in the plane $\{ z = 0 \}$ and not on axes (Figure \ref{rys:2}$(a)$),
		
		\item $P$ lies in the plane $\{ z = 0 \}$ and on the axis ${Ox}$
		(Figure \ref{rys:2}$(b)$).
	\end{enumerate}
\end{example}
\begin{remark}
	\label{rem5}Item \ref{it2} of Theorem \ref{Main} can be
	equivalently stated as follows:
	\begin{enumerate}
		\item[2'.] $\nu (\Gamma_+^P) < \nu (\Gamma_+)$ if, and only if, one of the
		following two conditions is satisfied:
		\begin{enumerate}
			\item $P$ lies in the interior of $\Gamma_-$ i.e.~$P \in \tmop{Int}
			(\Gamma_-)$,
			
			\item for each coordinate plane $H$ such that $P \in H$ the p-difference
			$\Gamma_+^P - \Gamma_+$ is either a pyramid with base $(\Gamma_+^P -
			\Gamma_+) \cap H$ and of height greater or equal to $2$, or an
			$n$-polytope with at least two vertices outside of $H$.
		\end{enumerate}
	\end{enumerate}
\end{remark}

\begin{example}
	The (weaker) requirement that the p-difference $\Gamma_+^P - \Gamma_+$ should lie
	in `a wall' of thickness $1$ around a coordinate plane is not sufficient for
	the equality $\nu (\Gamma_+^P) = \nu (\Gamma_+)$. In fact, if $\Gamma_+$ is
	the Newton polyhedron of the surface singularity $f (x, y, z) \assign x^6 +
	2 y^6 + z (x^2 + y^2) + z^4$ and $P = (3, 2, 0)$, then:
	\begin{enumerate}
		\item $\nu (\Gamma_+) = 15$, $\nu (\Gamma_+^P) = 13$,
		
		\item $\Gamma_+^P - \Gamma_+$ is a $3$-polytope with `base' $(\Gamma_+^P -
		\Gamma_+) \cap {Oxy}$ and of height $1$, but it is not a pyramid; it
		has two vertices above ${Oxy}$.
	\end{enumerate}
\end{example}

\section{Proof of the Theorem}

Let $\Gamma_+$ be a convenient Newton polyhedron in $\mathbbm{R}^3_{\geqslant
	0}$. Let $P \in
\mathbbm{N}_0^3 \setminus \Gamma_+$ denote a lattice point under $\Gamma_+$.

Item \ref{it1} of the theorem will be proved in the course of the proof of
item \ref{it2}, because we will in fact show that the negation of the
combinatorial condition in item \ref{it2} implies the strict inequality $\nu
(\Gamma_+^P) < \nu (\Gamma_+)$.

We first prove that the combinatorial condition in item \ref{it2} implies the equality $\nu(\Gamma_{+})=\nu(\Gamma_{+}^{P})$. Without loss of generality we may assume that, having fixed coordinates $(x, y, z)$ in $\mathbbm{R}^3$, we
have: $H = \{ z = 0 \}$, $P \in H$ and $\Gamma_+^P - \Gamma_+$ is a pyramid
with base $(\Gamma_+^P - \Gamma_+) \cap H$ and of height equal to $1$. We must
show that $\nu (\Gamma_+^P) = \nu (\Gamma_+)$. Since $\Gamma_+^P - \Gamma_+ =
\Gamma_- - \Gamma^P_-$, but the latter is a p-difference of two polytopes, we
prefer to reason in terms of $\Gamma_-$ and $\Gamma_-^P$ instead.

We have three possibilities:
\vspace{3mm}

$\textbf{1.}$ $\;$ $P$ does not lie on any axis, that is $P \nin {Ox}$ and $P \nin
{Oy}$. Then the polytopes $\Gamma_-$ and $\Gamma_-^P$ are identical on
${Ox}$, ${Oy}$, ${Oz}$, ${Oxz}$ and ${Oyz}$.
Denoting by $W$ the p-difference polygon of $\Gamma_-$ and $\Gamma^P_-$ in
${Oxy}$, we have by (\ref{for1})
\[ \nu (\Gamma_+) - \nu (\Gamma^P_+) = \frac{3! \cdummy \tmop{vol}_2 (W)
	\cdot 1}{3} - 2! \cdummy \tmop{vol}_2 (W) = 0. \]

\vspace{3mm}

$\textbf{2.}$ $\;$ $P$ lies on ${Ox}$ or ${Oy}$ and $P\neq 0$. Up to renaming of the
variables, we may assume that $P \in {Ox}$. Hence and by the assumption that $\Gamma_{+}$ is convenient the apex of the pyramid must lie in the plane $Oxz$. Denoting $W \assign
(\Gamma_- - \Gamma^P_-) \cap {Oxy}$ and $L \assign (\Gamma_- -
\Gamma^P_-) \cap {Ox}$, we have
\[ \nu (\Gamma_+) - \nu (\Gamma^P_+) = \frac{3! \cdummy \tmop{vol}_2 (W)
	\cdot 1}{3} - 2! \cdummy \tmop{vol}_2 (W) - \frac{2! \cdummy
	\tmop{vol}_1 (L) \cdot 1}{2} + \tmop{vol}_1 (L) = 0. \]

$\textbf{3.}$ $P=0$. Then $\Gamma_{-}^P=\varnothing$. By the assumption that $\Gamma_{+}$ is convenient the apex of the pyramid must be $(0,0,1)$. Hence, if we denote by $L_x$, $L_y$, $L_z$, $W_{xy}$, $W_{xz}$, $W_{yz}$ the intersections of $\Gamma_{-}$ with coordinate axes and plane, respectively, then
\begin{align*}
\nu(\Gamma_{+})-\nu(\Gamma_{+}^P)=\nu(\Gamma_{+})=\dfrac{3!\cdummy\tmop{vol}_2(W_{xy})\cdot 1}{3}-2!\cdummy \tmop{vol}_2(W_{xy})-2!\cdummy \tmop{vol}_2(W_{xz})-2!\cdummy \tmop{vol}_2(W_{yz})+\\
+\tmop{vol}_1(L_{x})+\tmop{vol}_1(L_{y})+\tmop{vol}_1(L_{z})-1=0.
\end{align*}

Let us pass to the proof of the inverse implication in item \ref{it2} of the
theorem. Assume to the contrary, that the combinatorial condition in item
\ref{it2} does not hold. Consider possible cases:

\vspace{3mm}

$\textbf{1.}$ $\;$ $P$ does not lie in any coordinate plane. Then the p-difference
$\Gamma_- - \Gamma^P_-$ is an $3$-polytope, disjoint from all the coordinate planes.
Hence,
\[ \nu (\Gamma_+) - \nu (\Gamma^P_+) = 3! \cdummy \tmop{vol}_3 (\Gamma_- -
\Gamma^P_-) > 0. \]

\vspace{3mm}

\begin{figure}[H]
	\centering
	\includegraphics{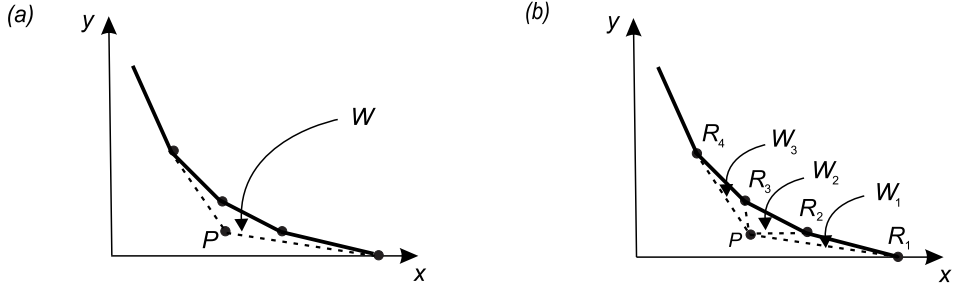}\\
	\caption{$(a)$ $p$-difference $\Gamma_{-}-\Gamma^{P}_{-}$ in $Oxy$. 
		$(b)$ Triangulation of $W$.}\label{rys:IIb}
\end{figure}

$\textbf{2.}$ \label{pII}$\;$ $P$ lies in a coordinate plane, but not on any axis.
Without loss of generality, we may assume that $P \in {Oxy} \backslash
({Ox} \cup {Oy})$. Then the p-difference $3$-polytope $\Gamma_- -
\Gamma^P_-$ is disjoint from the planes ${Oxz}$ and ${Oyz}$, but
$W \assign (\Gamma_- - \Gamma^P_-) \cap {Oxy} \neq \varnothing$ (Figure \ref{rys:IIb}$(a)$).
According to Remark \ref{rem5}, we should examine the following
possibilities:

\vspace{3mm}

$\textbf{a)}$ $\;$ $\Gamma_- - \Gamma^P_-$ is a pyramid with base $W$ and
of height $h \geqslant 2$. We have
\[ \nu (\Gamma_+) - \nu (\Gamma^P_+) = \frac{3! \cdummy \tmop{vol}_2 (W)
	\cdummy h}{3} - 2! \cdummy \tmop{vol}_2 (W) = 2 \cdummy \tmop{vol}_2 (W)
\cdummy (h - 1) > 0. \]

$\textbf{b)}$ \label{ppb} $\;$ There are at least two vertices of $\Gamma_-
- \Gamma^P_-$ lying above ${Oxy}$. Denote them by $Q_1, \ldots, Q_r$,
where $r \geqslant 2$. $W$ itself is a polygon in ${Oxy}$ of the shape
depicted in Figure  \ref{rys:IIb}$(a)$.

Naturally, $W$ is a union of relatively disjoint triangles $W_1, \ldots,
W_k$, $k \geqslant 1$, with one common vertex $P$. Let us enumerate all the
other vertices of $W$ according to their increasing $y$-coordinate and
denote them by $R_1, \ldots, R_{k + 1}$; then we may put $W_1 = \Delta (R_1,
R_2, P)$, $W_2 = \Delta (R_2, R_3, P)$, $\ldots$, $W_k = \Delta (R_k, R_{k +
	1}, P)$ (see Figure  \ref{rys:IIb}$(b)$).


We claim that each $W_i$ is the base of a pyramid $\mathbbm{P}_i$ with apex
in $\{ Q_1, \ldots, Q_r \}$ such that $\mathbbm{P}_i \subset (\Gamma_- -
\Gamma^P_-)$ and, moreover, all the $\mathbbm{P}_i$ are pairwise relatively
disjoint. In fact, each side $R_i R_{i + 1}$ is an edge of a 2-dimensional
facet of $\Gamma$, say $F_i$ $(i = 1, \ldots, k)$. Clearly, the facets $F_i$
are pairwise relatively disjoint and their vertices above ${Oxy}$ are
among $Q_1, \ldots, Q_r$. The pyramids $\mathbbm{P} (F_i, P)$ $(i = 1,
\ldots, k)$ are also pairwise relatively disjoint and lie in $\Gamma_- -
\Gamma_-^P$. Moreover, the pyramids are all convex, $F_i$ being convex.
For every $F_i$ choose arbitrarily $Q_{j_i}$ such that $Q_{j_i}$ is a
vertex of $F_i$ $(i = 1, \ldots, k)$. Define new pyramids $\mathbbm{P}_i \assign \mathbbm{P}
(W_i, Q_{j_i})$. Since $W_i$ is a facet of $\mathbbm{P} (F_i, P)$, we have
$\mathbbm{P}_i \subset \mathbbm{P} (F_i, P)$ and so $\mathbbm{P}_i$ are
pairwise relatively disjoint and lie in $\Gamma_- - \Gamma_-^P$ $(i = 1,
\ldots, k)$, as desired.

Let $h_i$ be the height of $\mathbbm{P}_i$. Clearly, $h_i \geqslant 1$ $(i =
1, \ldots, k)$. Moreover, the union $\bigcup_{1 \leqslant i \leqslant k}
\mathbbm{P}_i$ is not equal to $\Gamma_- - \Gamma_-^P$ because any
edge joining a pair of vertices from $Q_1, \ldots, Q_r$ (such edges always
exist) does not belong to any of the pyramids $\mathbbm{P}_1, \ldots,
\mathbbm{P}_k$, but such an edge is an edge of the $3$-polytope $\Gamma_- -
\Gamma_-^P$. Consequently, $V \assign (\Gamma_- - \Gamma_-^P) - \bigcup_{1
	\leqslant i \leqslant k} \mathbbm{P}_i$ is a non-empty compact $3$-polyhedron. We have 
\begin{align*}
\nu (\Gamma_+) - \nu (\Gamma^P_+) & =\nu (\Gamma_- - \Gamma^P_-)= 3! \cdummy \tmop{vol}_3 (V) + 3! \cdummy \sum_{1 \leqslant i
	\leqslant k} \frac{\tmop{vol}_2 (W_i) \cdummy h_i}{3} - 2! \cdummy \sum_{1
	\leqslant i \leqslant k} \tmop{vol}_2 (W_i)={}\\
& =3! \cdummy \tmop{vol}_3 (V) + 2! \cdummy \sum_{1 \leqslant i
	\leqslant k} \tmop{vol}_2 (W_i) \cdummy (h_i - 1) \geqslant 3! \cdummy \tmop{vol}_3
(V) > 0.
\end{align*}
\vspace{3mm}

$\textbf{3.}$ $\;$ $P$ lies on an axis and $P\neq 0$. For definiteness, let $P \in {Ox}$. Then the
p-difference $3$-polytope $\Gamma_- - \Gamma^P_-$ is disjoint from the plane
${Oyz}$ and from the axes ${Oy}$, ${Oz}$. Similarly as in $\textbf{2(b)}$
we divide the polygons $W \assign
(\Gamma_- - \Gamma^P_-) \cap {Oxy}$ and $\widetilde{W} \assign (\Gamma_- -
\Gamma^P_-) \cap {Oxz}$ into triangles $W_1, \ldots, W_k$, $k \geqslant
1$, and $\widetilde{W}_1, \ldots, \widetilde{W}_{\widetilde{k}}$, $\widetilde{k} \geqslant
1$, respectively, all of them having $P$ as the (only) common vertex.
\begin{figure}[H]
	\centering
	\includegraphics{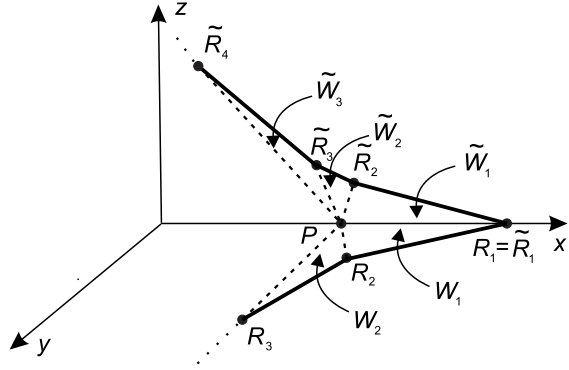}\\
	\caption{Triangulation of $W$ and $\widetilde{W}$.}\label{rys:III}
\end{figure}

We also arrange them such that $W_i = \Delta (R_i, R_{i + 1}, P)$ (resp.~$\widetilde{W}_j =
\Delta (\widetilde{R}_j, \widetilde{R}_{j + 1}, P)$), where the points $R_1, \ldots,
R_{k + 1}$ (resp.~$\widetilde{R}_1, \ldots, \widetilde{R}_{\widetilde{k} + 1}$) are
enumerated according to their increasing $y$- (resp. $z$-) coordinates (see
Figure \ref{rys:III}). Observe that, in particular, $R_1 = \widetilde{R}_1$.

\begin{claim}
	$W_1, \ldots,
	W_k$ and $\widetilde{W}_1, \ldots, \widetilde{W}_{\widetilde{k}}$ are bases of some
	$3$-pyramids $\mathbbm{P}_1, \ldots, \mathbbm{P}_k$ and
	$\widetilde{\mathbbm{P}}_1, \ldots, \widetilde{\mathbbm{P}}_{\widetilde{k}}$,
	respectively, all of which are pairwise relatively disjoint, possibly except
	for the pair $\mathbbm{P}_1$ and $\widetilde{\mathbbm{P}}_1$, lie in $\Gamma_- -
	\Gamma_-^P$ and whose vertices are taken from the set of vertices of
	$\Gamma_+$. Moreover, if $\mathbbm{P}_1$ and $\widetilde{\mathbbm{P}}_1$ are not
	relatively disjoint, then the triangle $\Delta (R_1, R_2, \widetilde{R}_2)$ is a facet of
	$\Gamma$ and $\mathbbm{P}_1 = \widetilde{\mathbbm{P}}_1 = \Delta (R_1, R_2,
	\widetilde{R}_2, P)$.
\end{claim}

\begin{claimproof}
	The segments $R_1 R_2$, $R_2 R_3$,
	$\ldots$ are edges of uniquely determined facets, say $F_1, \ldots,\allowbreak F_k$, of
	$\Gamma$. Clearly, $F_1, \ldots, F_k$ are pairwise relatively disjoint.
	Similarly, $\widetilde{R}_1 \widetilde{R}_2$, $\widetilde{R}_2 \widetilde{R}_3$, $\ldots$
	are edges of uniquely determined facets, say $\widetilde{F}_1, \ldots,
	\widetilde{F}_{\widetilde{k}}$, of $\Gamma$. It may happen that some $F_i$ are equal
	to some $\widetilde{F}_j$ (see Figure  \ref{rys:claim1}$(a)$).
	\begin{figure}[H]
		\centering
		\includegraphics{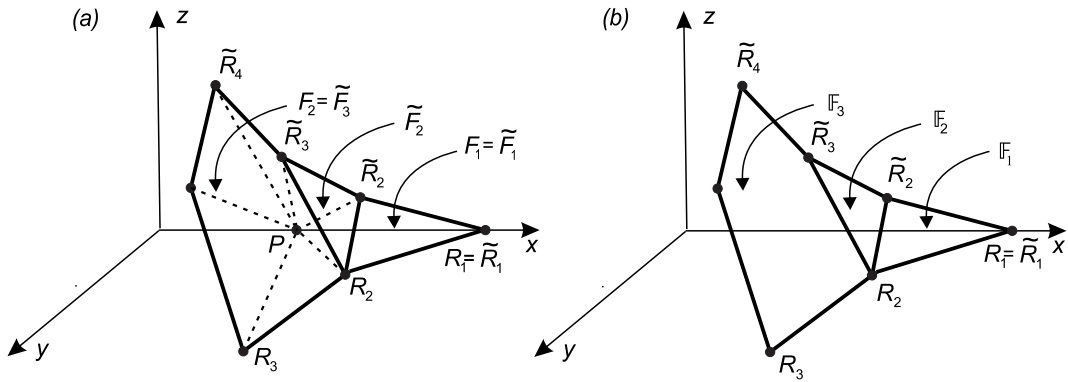}\\
		\caption{Facets of $\Gamma$ in $\Gamma_{-}-\Gamma^{P}_{-}$.}\label{rys:claim1}
	\end{figure}  
	Removing the duplicated facets from the sequence $F_1, \ldots, F_k,
	\widetilde{F}_1, \ldots, \widetilde{F}_{\widetilde{k}}$, we obtain a new sequence
	$\mathbbm{F}_1, \ldots, \mathbbm{F}_l$ of pairwise relatively disjoint
	facets of $\Gamma$ (Figure  \ref{rys:claim1}$(b)$) having the following properties:
	\begin{itemize}
		\item each side $R_i R_{i + 1}$ $(i = 1, \ldots, k)$ is an edge of a unique
		$\mathbbm{F}_s$ and, similarly, each side $\widetilde{R}_j
		\widetilde{R}_{j + 1}$ $(j = 1, \ldots, \widetilde{k})$ is an edge of a unique
		$\mathbbm{F}_{\widetilde{s}}$,\  $s,\widetilde{s}\in\{1,\ldots,l\}$,
		
		\item each $\mathbbm{F}_s$ $(s = 1, \ldots, l)$ has either one or at most
		two of its edges among $R_1 R_2$, $\ldots, R_k R_{k + 1}$, $\widetilde{R}_1
		\widetilde{R}_2$, $\ldots$, $\widetilde{R}_{\widetilde{k}} \widetilde{R}_{\widetilde{k} + 1}$.
	\end{itemize}
	Upon renaming the $\mathbbm{F}_s$, we may assume that $\mathbbm{F}_1$ has
	$R_1 R_2$ as one of its edges. For each $s \in \{ 1, \ldots, l \}$ we build
	the $3$-pyramid $\mathbbm{P} (\mathbbm{F}_s, P)$. These pyramids are
	convex, pairwise relatively disjoint, lie in $\Gamma_- - \Gamma_-^P$ and each of
	them has either one or at most two of its edges among $R_1 R_2$, $\ldots,
	R_k R_{k + 1}$, $\widetilde{R}_1 \widetilde{R}_2$, $\ldots$, $\widetilde{R}_{\widetilde{k}}
	\widetilde{R}_{\widetilde{k} + 1}$. Hence, each $\mathbbm{P} (\mathbbm{F}_s, P)$ $(s
	= 1, \ldots, l)$ has either one or at most two of its facets among $W_1$,
	$\ldots, W_k$, $\widetilde{W}_1$, $\ldots$, $\widetilde{W}_{\widetilde{k}}$. In order to
	prove the claim, we only need to build the $\mathbbm{P}_1, \ldots,
	\mathbbm{P}_k$, $\widetilde{\mathbbm{P}}_1, \ldots,
	\widetilde{\mathbbm{P}}_{\widetilde{k}}$ as appropriately chosen pyramids hiding
	inside the larger $\mathbbm{P} (\mathbbm{F}_s, P)$ $(s = 1, \ldots, l)$.
	
	Fix $s_0\in\{1,\ldots,l\}$.  If $\mathbbm{P} (\mathbbm{F}_{s_0}, P)$ has only one facet among $W_1$,
	$\ldots, W_k$, $\widetilde{W}_1$, $\ldots$, $\widetilde{W}_{\widetilde{k}}$, then we
	simply consider a pyramid whose base is equal to this facet and whose apex
	is chosen as any other vertex of $\mathbbm{F}_{s_0}$; such pyramid is then
	contained in $\mathbbm{P} (\mathbbm{F}_{s_0}, P)$.
	
	If $\mathbbm{P} (\mathbbm{F}_{s_0}, P)$ has two facets among $W_1$, $\ldots,
	W_k$, $\widetilde{W}_1$, $\ldots$, $\widetilde{W}_{\widetilde{k}}$, then, by our
	construction, one of them is equal to some $W_p$, where $p \in \{ 1, \ldots,
	k \}$, and the other one is equal to some $\widetilde{W}_q$, where $q \in \{ 1,
	\ldots, \widetilde{k} \}$. Excluding $s_0 = 1$, we notice that there exist two
	relatively disjoint pyramids contained in $\mathbbm{P} (\mathbbm{F}_{s_0},
	P)$: one has its base equal to $W_p$, the other one has base $\widetilde{W}_q$,
	and both of them have their apices appropriately chosen from
	$\mathbbm{F}_{s_0}$. This follows from the fact that if $s_0 \geqslant 2$,
	then the facet $\mathbbm{F}_{s_0}$ has at least four vertices (see Figure \ref{rys:claim1}$(b)$),
	$R_p$, $R_{p + 1}$, $\widetilde{R}_q, \widetilde{R}_{q + 1}$, which allows us to construct relatively disjoint pyramids $\mathbbm{P}_{p}\assign\Delta (P,
	R_p, R_{p + 1}, \widetilde{R}_q)$ of base $W_{p}$ and $\widetilde{\mathbbm{P}}_{q}\assign\Delta (P, \widetilde{R}_q, \widetilde{R}_{q +
		1}, R_{p+1})$ of base $\widetilde{W}_{q}$, both of them contained inside the pyramid $\mathbbm{P}
	(\mathbbm{F}_{s_0}, P)$ (see Figure \ref{rys:claim2}).
	\begin{figure}[H]
		\centering
		\includegraphics{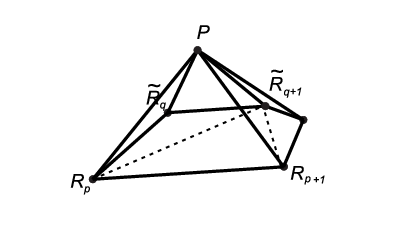}\\
		\caption{The pyramid $\mathbb{P}(\mathbb{F}_{s_{0}},P)$ turned upside down with inscribed two pyramids of bases $W_{p}$ and $\widetilde{W}_{q}$.}\label{rys:claim2}
	\end{figure}
	
	Now, let us treat the case $s_0 = 1$ and $\mathbbm{F}_1$ has two of its
	facets among $W_1$, $\ldots, W_k$, $\widetilde{W}_1$, $\ldots$,
	$\widetilde{W}_{\widetilde{k}}$. According to our earlier arrangements, one of these
	facets is equal to $W_1$ and the other one has to be equal to $\widetilde{W}_1$.
	If $\mathbbm{F}_1$ is not equal to the triangle $\Delta (R_1, R_2,
	\widetilde{R}_2)$, i.e. $\mathbbm{F}_1$ has a fourth vertex, say $Q$, we may
	take $\mathbbm{P}_1 \assign \mathbbm{P} (W_1, Q)$, $\widetilde{\mathbbm{P}}_1
	\assign \mathbbm{P} (\widetilde{W}_1, Q)$. These pyramids are relatively
	disjoint and lie in $\mathbbm{P} (\mathbbm{F}_1, P)$. Such choice of
	pyramids is impossible if $\mathbbm{F}_1 = \Delta (R_1, R_2, \widetilde{R}_2)$,
	hence in this case we simply put $\mathbbm{P}_1 = \widetilde{\mathbbm{P}}_1
	\assign \Delta (R_1, R_2, \widetilde{R}_2, P) 
	=\mathbbm{P} (W_1, \widetilde{R}_2) =\mathbbm{P} (\widetilde{W}_1, R_2)$.
	
	Since each $W_1, \ldots, W_k$, $\widetilde{W}_1, \ldots, \widetilde{W}_{\widetilde{k}}$
	is a facet of some $\mathbbm{P} (\mathbbm{F}_s, P)$, where $s \in \{ 1,
	\ldots, l \}$, we have constructed pyramids $\mathbbm{P}_{1},\ldots,\mathbbm{P}_{k},\widetilde{\mathbbm{P}}_{1},\ldots,\widetilde{\mathbbm{P}}_{\widetilde{k}}$ with all the required
	properties.
\end{claimproof}

Denote by $h_i$ and $\widetilde{h}_j$ the heights of $\mathbbm{P}_i$ and
$\widetilde{\mathbbm{P}}_j$, respectively. They are positive integers. Take the $p$-difference $V
\assign (\Gamma_- - \Gamma_-^P) - \bigcup_{1 \leqslant i \leqslant k}
\mathbbm{P}_i - \bigcup_{1 \leqslant j \leqslant \widetilde{k}}
\widetilde{\mathbbm{P}}_j$. Then $V$ is either empty or it is a compact $3$-polyhedron. Let $L
\assign R_1 P = (\Gamma_- - \Gamma_-^P) \cap {Ox}$, that is $L$ be the
p-difference $(\Gamma_- - \Gamma_-^P)$ on the axis ${Ox}$.

According to the above claim, we have two possibilities to consider:
\vspace{3mm}

$\textbf{a)}$ $\;$ $\mathbbm{P}_1 \neq \widetilde{\mathbbm{P}}_1$. Using
formula (\ref{for1}), we compute
\begin{align*}
\nu (\Gamma_+) - \nu (\Gamma^P_+)  &=\!\begin{multlined}[t][.7\displaywidth]
3! \cdummy \tmop{vol}_3 (V) + 3!
\cdummy \sum_{1 \leqslant i \leqslant k} \tmop{vol}_3 (\mathbbm{P}_i) + 3!
\cdummy \sum_{1 \leqslant j \leqslant \widetilde{k}} \tmop{vol}_3
(\widetilde{\mathbbm{P}}_j)+{}\\     - 2! \cdummy \sum_{1 \leqslant i \leqslant k} \tmop{vol}_2 (W_i) -
2! \cdummy \sum_{1 \leqslant j \leqslant \widetilde{k}} \tmop{vol}_2
(\widetilde{W}_j) + 1! \cdummy \tmop{vol}_1 (L)=
\end{multlined}\\
& =\!\begin{multlined}[t][.7\displaywidth]
3! \cdummy \tmop{vol}_3 (V) + 2 \cdummy \sum_{1 \leqslant i
	\leqslant k} \tmop{vol}_2 (W_i) \cdummy (h_i - 1)+{}\\
+ 2 \cdummy \sum_{1 \leqslant j \leqslant \widetilde{k}} \tmop{vol}_2
(\widetilde{W}_j) \cdummy (\widetilde{h}_j - 1) + \tmop{vol}_1 (L)
\geqslant  \tmop{vol}_1 (L) > 0.
\end{multlined}
\end{align*}
This gives the required inequality.  Notice that in this case we actually do not use the assumption that
$\Gamma_- - \Gamma_-^P$ is not a pyramid of height one.

\vspace{3mm}

$\textbf{b)}$ $\;$ $\mathbbm{P}_1 = \widetilde{\mathbbm{P}}_1 = \Delta (R_1,
R_2, \widetilde{R}_2, P)$ and $\Delta (R_1, R_2, \widetilde{R}_2)$ is a face of
$\Gamma$. We have
\begin{align*}
\MoveEqLeft \nu(\Gamma_+) - \nu (\Gamma^P_+) =\\
& =\!\begin{multlined}[t][.86\displaywidth]
3! \cdummy \tmop{vol}_3 (V) + 3!
\cdummy \sum_{1 \leqslant i \leqslant k} \tmop{vol}_3 (\mathbbm{P}_i) + 3!
\cdummy \sum_{2 \leqslant j \leqslant \widetilde{k}} \tmop{vol}_3
(\widetilde{\mathbbm{P}}_j)+{}\\
- 2! \cdummy \sum_{1 \leqslant i \leqslant k} \tmop{vol}_2 (W_i) -
2! \cdummy \sum_{1 \leqslant j \leqslant \widetilde{k}} \tmop{vol}_2
(\widetilde{W}_j) + 1! \cdummy \tmop{vol}_1 (L)=
\end{multlined}\\
& = \!\begin{multlined}[t][.86\displaywidth]
3! \cdummy \tmop{vol}_3 (V) + 2 \cdummy \sum_{2 \leqslant i
	\leqslant k} \tmop{vol}_2 (W_i) \cdummy (h_i - 1) + 2 \cdummy \sum_{2
	\leqslant j \leqslant \widetilde{k}} \tmop{vol}_2 (\widetilde{W}_j) \cdummy
(\widetilde{h}_j - 1)+{}\\
+ 2 \cdummy \tmop{vol}_2 (W_1) \cdummy (h_1 - 1) - 2 \cdummy \tmop{vol}_2
(\widetilde{W}_1) + \tmop{vol}_1 (L) .
\end{multlined}
\end{align*}
Since $W_1 =\mathbbm{P} (L, R_2)$, $\widetilde{W}_1 =\mathbbm{P} (L,
\widetilde{R}_2)$ are two perpendicular facets of $\Delta (R_1, R_2,
\widetilde{R}_2, P)$, we can continue the above equality

\begin{align*}
& =\!\begin{multlined}[t][.86\displaywidth]
3! \cdummy \tmop{vol}_3 (V) + 2 \cdummy \sum_{2 \leqslant i
	\leqslant k} \tmop{vol}_2 (W_i) \cdummy (h_i - 1) + 2 \cdummy \sum_{2
	\leqslant j \leqslant \widetilde{k}} \tmop{vol}_2 (\widetilde{W}_j) \cdummy
(\widetilde{h}_j - 1){}
+ 2 \cdummy \frac{\tmop{vol}_1 (L)}{2} \cdummy \widetilde{h}_1
\cdummy (h_1 - 1) + \\- 2 \cdummy \frac{\tmop{vol}_1 (L)}{2} \cdummy h_1 +
\tmop{vol}_1 (L) = \\
\end{multlined}\\
& = \!\begin{multlined}[t][.86\displaywidth]
3! \cdummy \tmop{vol}_3 (V) + 2 \cdummy \sum_{2 \leqslant i
	\leqslant k} \tmop{vol}_2 (W_i) \cdummy (h_i - 1) + 2 \cdummy \sum_{2
	\leqslant j \leqslant \widetilde{k}} \tmop{vol}_2 (\widetilde{W}_j) \cdummy
(\widetilde{h}_j - 1)+{}
\tmop{vol}_1 (L) \cdummy (\widetilde{h}_1 - 1) \cdummy (h_1 - 1). 
\end{multlined}
\end{align*}


This shows that, whatever the situation, it holds $\nu (\Gamma_+) - \nu
(\Gamma^P_+) \geqslant 0$ in this case. Moreover, $\nu (\Gamma_+) - \nu
(\Gamma^P_+) > 0$ if, and only if, one of the following three conditions is
satisfied:
\begin{enumerateroman}
	\item \label{red1}$V \neq \varnothing$,
	
	\item \label{red2} there exists $h_i$, $i\in\{2,\ldots,k\}$ or $\widetilde{h}_j$, $j\in\{2,\ldots,\widetilde{k}\}$ greater or equal to $2$, 
	
	\item \label{red3}both $h_1$ and $\widetilde{h}_1$ are greater or equal to
	$2$.
\end{enumerateroman}
Hence, the proof of the theorem will be finished once we show that the
assumption `$\Gamma_- - \Gamma_-^P$ is not a pyramid of height one' implies
at least one of these conditions.

\begin{figure}[H]
	\centering
	\includegraphics{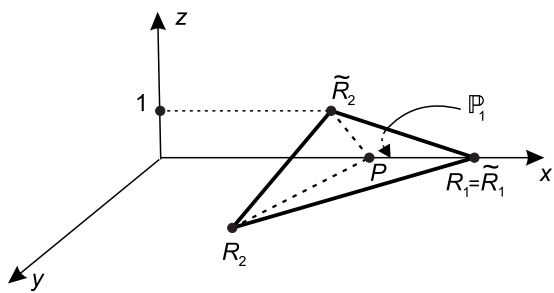}\\
	\caption{The pyramid $\mathbbm{P}_1$ of height $1$.}\label{rys:IIIb}
\end{figure}

To this end assume that (\ref{red3}) is not satisfied i.e. $h_1 = 1$ or
$\widetilde{h}_1 = 1$, say $h_1 = 1$ (see Figure \ref{rys:IIIb}).
Since $\mathbbm{P}_1 = \mathbbm{P}(W_1,\widetilde{R}_2)$ is a pyramid with
base $W_1$ contained in a coordinate plane and of height $1$, our assumption
implies that $\Gamma_- - \Gamma^P_- \neq \mathbbm{P}_1$. Hence, there exists
a facet $F \in \Gamma$ of $\Gamma_- - \Gamma^P_-$ sharing the edge
$R_2 \widetilde{R}_2$ in common with the facet $\mathbbm{F}_1 = \Delta (R_1, R_2, \widetilde{R}_2)_{}$.

\begin{figure}[H]
	\centering
	\includegraphics{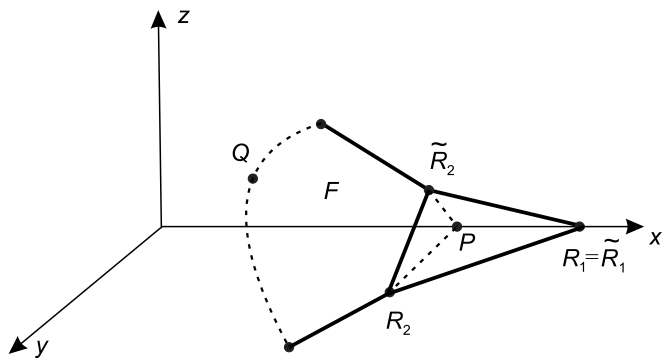}\\
	\caption{The facet $F$ with a vertex $Q$ outside $Oxy\cup Oxz$}\label{rys:IIIb2}
\end{figure}

If $F$ has some vertex $Q$ outside of ${Oxy} \cup {Oxz}$ (see Figure \ref{rys:IIIb2}), then
(\ref{red1}) holds because $F$ is not entirely covered by the pyramids $\mathbbm{P}_2,\ldots, \mathbbm{P}_k, \widetilde{\mathbbm{P}}_2,\ldots, \widetilde{\mathbbm{P}}_{\widetilde{k}}$ (see explanations below) which implies $F-\bigcup_{2
	\leqslant i \leqslant k}\mathbbm{P}_i-\bigcup_{2 \leqslant j \leqslant
	\widetilde{k}}\widetilde{\mathbbm{P}}_j$ lies in the boundary of $(\Gamma_- - \Gamma_-^P) - \bigcup_{1
	\leqslant i \leqslant k} \mathbbm{P}_i - \bigcup_{1 \leqslant j \leqslant
	\widetilde{k}} \widetilde{\mathbbm{P}}_j$ i.e. in the boundary of $V$. Hence $V \neq \varnothing$ i.e.~(\ref{red1}) holds.
To see $F$ is not entirely covered by the pyramids $\mathbbm{P}_2,\ldots, \mathbbm{P}_k, \widetilde{\mathbbm{P}}_2,\ldots, \widetilde{\mathbbm{P}}_{\widetilde{k}}$ we consider the three possible cases:

$\textbf{1.}$ $F$ has no edges in ${Oxy} \cup {Oxz}$. Then $F$ is not equal to any $F_2,\ldots,F_k,\widetilde{F}_2,\ldots,\widetilde{F}_{\widetilde{k}}$ considered it the proof of Claim. Hence by construction of $\mathbbm{P}_2,\ldots, \mathbbm{P}_k, \widetilde{\mathbbm{P}}_2,\ldots, \widetilde{\mathbbm{P}}_{\widetilde{k}}$ the facet $F$ is not covered by them. Precisely $F-\bigcup_{2
	\leqslant i \leqslant k}\mathbbm{P}_i-\bigcup_{2 \leqslant j \leqslant
	\widetilde{k}}\widetilde{\mathbbm{P}}_j=F$.

$\textbf{2.}$ $F$ has only one edge in ${Oxy} \cup {Oxz}$. Then $F$ is equal either to $F_2$ or to $\widetilde{F}_2$.  Since in this case $F$ has at least $4$ vertices the pyramid either $\mathbbm{P}_{2}$ or $\widetilde{\mathbbm{P}}_{2}$, inscribed in $\mathbbm{P}(F,P)$ according to our construction of $\mathbbm{P}_1,\ldots, \mathbbm{P}_k, \widetilde{\mathbbm{P}}_1,\ldots, \widetilde{\mathbbm{P}}_{\widetilde{k}}$, does not cover entirely $F$ (notice $\mathbbm{P}_2$ and $\widetilde{\mathbbm{P}}_2$ has only triangles as facets). Consequently $F-\mathbbm{P}_{2}$ or $F-\widetilde{\mathbbm{P}}_{2}$ is not empty.

$\textbf{3.}$ $F$ has two edges in ${Oxy} \cup {Oxz}$, necessarily one in ${Oxy}$ and one in ${Oxz}$. Then $F=\mathbbm{F}_i$ for some $i\in\{2,\ldots,l\}$. Since in this case $F$ has at least $5$ vertices and two pyramids $\mathbbm{P}_{i'}$ and $\widetilde{\mathbbm{P}}_{j}$ (for some $i'\in\{2,\ldots,k\}$ and $j\in\{2,\ldots,\widetilde{k}\}$), inscribed in $\mathbbm{P}(\mathbbm{F}_i,P)$, have triangles as facets then $F-\mathbbm{P}_{i'}-\widetilde{\mathbbm{P}}_{j}$ is not empty.

If all the vertices of $F$ belong to ${Oxy} \cup {Oxz}$, then
there are only three possible scenarios:
\begin{itemize}
	\item[$(\star)$] $F$ is a convex quadrilateral; then apart from $R_2$ and
	$\widetilde{R}_2$, it has only two further vertices, necessarily $R_3$ and
	$\widetilde{R}_3$ (see Figure \ref{rys:IIIb3}$ (a) $).
	\begin{figure}[H]
		\centering
		\includegraphics{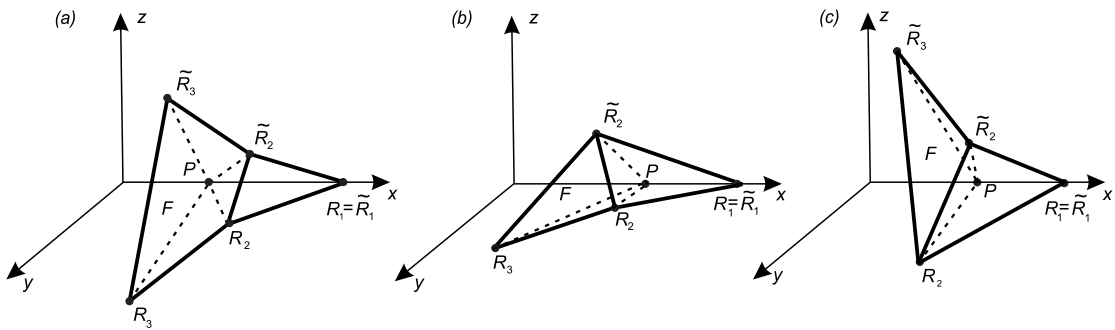}\\
		\caption{The possible placements of the face $F$ with vertices in $Oxy\cup Oxz$. }\label{rys:IIIb3}
	\end{figure}    
	\item[$(\star\star)$] $F$ is a triangle with the third vertex lying in
	${Oxy}$; this vertex is necessarily $R_3$ (see Figure \ref{rys:IIIb3}(b)).
	\item[$({\star\star}\star)$] $F$ is a triangle with the third vertex lying
	in ${Oxz}$; this vertex is necessarily $\widetilde{R}_3$ (see Figure \ref{rys:IIIb3}(c)).
\end{itemize}
In case ($\star$) our construction of $\mathbbm{P}_i$,
$\widetilde{\mathbbm{P}}_j$ shows that either $\mathbbm{P}_2$ or
$\widetilde{\mathbbm{P}}_2$ has its height greater than or equal to $2$. Hence,
(\ref{red2}) is satisfied.

In case ($\star\star$) we are in a position to repeat the above
construction of $F$, this time for the edge $R_3 \widetilde{R}_2$ instead of
$R_2 \widetilde{R}_2$, which leads to some new facet $F'$. To this facet $F'$
also applies the same case analysis as the one previously performed for
$F$.

In case (${\star\star}\star$) we have two options: the first one $\widetilde{h}_1 \geqslant
2$; then our construction of $\widetilde{\mathbbm{P}}_2$ shows that we must have
$\widetilde{h}_2 = \widetilde{h}_1 \geqslant 2$ so we get condition (\ref{red2}).
If, however, $\widetilde{h}_1 = 1$, we are in the position to repeat the above
construction of $F$, this time for the edge $R_2 \widetilde{R}_3$ instead of
$R_2 \widetilde{R}_2$, which leads to some new facet $F'$. To this facet $F'$
also applies the same case analysis as the one previously performed for $F$.

We see that, after finitely many steps, we either find a facet $F \in
\Gamma$ of $\Gamma_- - \Gamma^P_-$ with one of its vertices outside of
${Oxy} \cup {Oxz}$ (then $V \neq \varnothing$ and then (\ref{red1})
holds), or we discover that (\ref{red2}) holds. 

$\textbf{4.}$ $P=0$. Then $\Gamma_{-}^P=\varnothing$ and consequently $\nu(\Gamma_{+})-\nu(\Gamma_{+}^P)=\nu(\Gamma_{-})-\nu(\Gamma_{-}^P)=\nu(\Gamma_{-})$. By assumption $\Gamma_{-}$ is not a pyramid of base in a coordinate plane and height $1$. Hence $\Gamma_{+}$ intersects axes at points with coordinates greater or equal to $2$. Take auxiliary point $\widetilde{P}=(1,0,0)$ on axis $Ox$ liying under $\Gamma_{+}$. Then $\Gamma_{-}-\Gamma_{-}^{\widetilde{P}}$ is also not a pyramid of height $1$. Hence by \textbf{3} $\nu(\Gamma_{+})-\nu(\Gamma_{+}^{\widetilde{P}})>0$. Since $\Gamma_{-}^{\widetilde{P}}$ is a pyramid of height $1$ then $\nu(\Gamma_{-}^{\widetilde{P}})=0$. But $\Gamma_{-}=(\Gamma_{-}-\Gamma_{-}^{\widetilde{P}})\cup \Gamma_{-}^{\widetilde{P}}$ and of course $\nu((\Gamma_{-}-\Gamma_{-}^{\widetilde{P}})\cup \Gamma_{-}^{\widetilde{P}})=\nu(\Gamma_{-}-\Gamma_{-}^{\widetilde{P}})+\nu(\Gamma_{-}^{\widetilde{P}})$. Hence $\nu(\Gamma_{-})=\nu(\Gamma_{-}-\Gamma_{-}^{\widetilde{P}})=\nu(\Gamma_{-})-\nu(\Gamma_{-}^{\widetilde{P}})=\nu(\Gamma_{+})-\nu(\Gamma_{+}^{\widehat{P}})>0$. This ends the proof.

\begin{corollary}
	Let $\Gamma_{+}$ be a convenient Newton polyhedron in $\mathbb{R}^3_{\geqslant 0}$. Then $\nu(\Gamma_{+})\geqslant 0$. Moreover, $\nu(\Gamma_{+})=0$ if and only if $\Gamma_{+}=\mathbb{R}^3_{\geqslant 0}$ or $\Gamma_{+}$ intersects one of axes at the point with coordinate equal to $1$.
\end{corollary}

\begin{remark}
	The last corollary was proved by M. Furuya \cite{Fur04}, Corollary $2.4$, in $n$-dimensional case.
\end{remark}


\bibliographystyle{abbrv}
\bibliography{bibliografia}   

\begin{thebibliography}{10}

\bibitem{Arn04}
V.~I. Arnold.
\newblock {\em Arnold's problems}.
\newblock Springer-Verlag, Berlin; PHASIS, Moscow, 2004.
\newblock Translated and revised edition of the 2000 Russian original, With a
  preface by V. Philippov, A. Yakivchik and M. Peters.

\bibitem{Ber09}
M.~Berger.
\newblock {\em Geometry {II}}.
\newblock Universitext. Springer-Verlag, Berlin, 1987.
\newblock Translated from the French by M. Cole and S. Levy.

\bibitem{Biv09}
C.~Bivi\`a-Ausina.
\newblock Local \l ojasiewicz exponents, {M}ilnor numbers and mixed
  multiplicities of ideals.
\newblock {\em Math. Z.}, 262(2):389--409, 2009.

\bibitem{Fur04}
M.~Furuya.
\newblock Lower bound of {N}ewton number.
\newblock {\em Tokyo J. Math.}, 27(1):177--186, 2004.

\bibitem{GN12}
G.-M. Greuel and H.~D. Nguyen.
\newblock Some remarks on the planar {K}ouchnirenko's theorem.
\newblock {\em Rev. Mat. Complut.}, 25(2):557--579, 2012.

\bibitem{Gwo08}
J.~Gwo\'zdziewicz.
\newblock Note on the {N}ewton number.
\newblock {\em Univ. Iagel. Acta Math.}, (46):31--33, 2008.

\bibitem{Kou76}
A.~G. Kouchnirenko.
\newblock Poly\`edres de {N}ewton et nombres de {M}ilnor.
\newblock {\em Invent. Math.}, 32(1):1--31, 1976.

\bibitem{Len08}
A.~Lenarcik.
\newblock On the {J}acobian {N}ewton polygon of plane curve singularities.
\newblock {\em Manuscripta Math.}, 125(3):309--324, 2008.

\bibitem{Mil68}
J.~Milnor.
\newblock {\em Singular points of complex hypersurfaces}.
\newblock Annals of Mathematics Studies, No. 61. Princeton University Press,
  Princeton, N.J., 1968.

\bibitem{Mon16}
P.~Mondal.
\newblock {Intersection multiplicity, Milnor number and Bernstein's theorem}.
\newblock {\em ArXiv}, 2016.

\bibitem{Oka89}
M.~Oka.
\newblock On the weak simultaneous resolution of a negligible truncation of the
  {N}ewton boundary.
\newblock In {\em Singularities ({I}owa {C}ity, {IA}, 1986)}, volume~90 of {\em
  Contemp. Math.}, pages 199--210. Amer. Math. Soc., Providence, RI, 1989.

\bibitem{Ste85}
J.~H.~M. Steenbrink.
\newblock Semicontinuity of the singularity spectrum.
\newblock {\em Invent. Math.}, 79(3):557--565, 1985.

\bibitem{Var84}
A.~N. Varchenko.
\newblock Asymptotic behaviors of integrals, and {H}odge structures.
\newblock In {\em Current problems in mathematics, {V}ol. 22}, Itogi Nauki i
  Tekhniki, pages 130--166. Akad. Nauk SSSR, Vsesoyuz. Inst. Nauchn. i Tekhn.
  Inform., Moscow, 1983.

\end{thebibliography}


\end{document}